\documentclass[a4paper,reqno]{amsart}

\usepackage[%
backend=bibtex,%
style=alphabetic,%
natbib=false,%
isbn=false, 
doi=false, 
url=false, 
eprint=true, 
giveninits=true, 
backref=false, 
maxbibnames=99,%
]{biblatex}%

\AtEveryBibitem{%
  \ifentrytype{unpublished}{}{\clearfield{note}}%
  }%
\renewbibmacro{in:}{} 

\usepackage{amssymb,amsmath}%
\usepackage{doi}%
\usepackage{eucal}%
\usepackage{graphicx}%
\usepackage{mathscinet}%
\usepackage{mathtools}%
\usepackage{tikz-cd}%
\usetikzlibrary{patterns}
\usepackage{xparse}%

\usepackage{hyperref}%
\hypersetup{%
  colorlinks,%
  linkcolor={red!80!black},%
  citecolor={blue!80!black},%
  urlcolor={blue!80!black}%
}%

\usepackage{thmtools}%
\usepackage{cleveref}%

\numberwithin{equation}{section}%

\newcounter{thmintro}%

\declaretheorem[style=theorem,sibling=equation]{theorem}%
\declaretheorem[style=theorem,name=Theorem]{theorem*}%
\declaretheorem[style=definition,numbered=no]{acknowledgements}%
\declaretheorem[style=definition,numbered=no,name={Financial Support}]{financialsupport}%
\declaretheorem[style=remark,sibling=equation]{remark}%
\declaretheorem[style=remark,sibling=equation]{example}%

\newcommand{\kk}{\mathbf{k}}%
\newcommand{\ZZ}{\mathbb{Z}}%
\DeclareMathOperator{\chark}{char}

\NewDocumentCommand{\id}{sO{}}{\IfBooleanTF{#1}{1}{\mathrm{id}}_{#2}}%

\NewDocumentCommand{\mmod}{s m}{\operatorname{mod}\IfBooleanTF{#1}{(#2)}{#2}} 

\NewDocumentCommand{\Hom}{O{}mm}{\operatorname{Hom}_{#1}(#2,#3)}%
\NewDocumentCommand{\RHom}{sO{}D<>{\bullet,*}mm}{\mathbb{R}\!\operatorname{Hom}_{#2}\IfBooleanT{#1}{^{#3}}(#4,#5)}%
\NewDocumentCommand{\dgHom}{O{}mmO{}}{\mathbf{Hom}_{#1}^{#4}(#2,#3)}%
\NewDocumentCommand{\Ext}{O{}O{\bullet}mm}{\operatorname{Ext}_{#1}^{#2}(#3,#4)}

\NewDocumentCommand{\DerCat}{O{}mD<>{}}{\operatorname{D}\IfValueT{#1}{^{\mathrm{#1}}}_{\mathrm{#3}}(#2)}%

\newcommand{\category}[1]{\mathcal{#1}}%
\renewcommand{\O}{\category{O}}%
\renewcommand{\P}{\category{P}}%
\newcommand{\T}{\category{T}}%

\RenewDocumentCommand{\H}{O{\bullet}m}{\operatorname{H}^{#1}(#2)}%
\NewDocumentCommand{\Hclass}{m}{\{ \!\!\{ {#1} \} \!\! \}}

\NewDocumentCommand{\HC}{sO{\bullet}D<>{*}mod<>}{\operatorname{C}\IfBooleanTF{#1}{^{#2}}{^{#2,#3}}\!\left(#4\IfValueT{#5}{,#5}\right)\IfValueT{#6}{\![#6]}}
\NewDocumentCommand{\HH}{sO{\bullet}D<>{*}mod<>}{\operatorname{HH}\IfBooleanTF{#1}{^{#2}}{^{#2,#3}}\!\left(#4\IfValueT{#5}{,#5}\right)\IfValueT{#6}{\![#6]}}

\newcommand{\Hd}{d_{\mathrm{Hoch}}}

\NewDocumentCommand{\cupp}{mm}{{#1}\cdot{#2}}
\NewDocumentCommand{\braces}{mm}{{#1}\left\{#2\right\}}
\DeclareMathOperator{\Sq}{Sq}

\NewDocumentCommand{\Astr}{sD<>{}O{A}}{\IfBooleanTF{#1}{\overline{m}}{m}_{#2}^{#3}}

\NewDocumentCommand{\Map}{O{}mm}{\operatorname{Map}_{#1}\!\left(#2,#3\right)}

\NewDocumentCommand{\opA}{O{\infty}}{\mathbb{A}_{#1}}
\NewDocumentCommand{\opEnd}{m}{\mathcal{E}\!\left(#1\right)}
\NewDocumentCommand{\BK}{O{r}O{\bullet}D<>{*}}{E_{#1}^{#2,#3}}

\newcommand{\bfLambda}{\mathbf{\Lambda}}

\newcounter{rpage}
\newcommand{\rpage}{\value{rpage}} 				
\newcommand{\rescale}{0.3}			
\newcommand{\margen}{0.1}			
\newcommand{\abajo}{5}				
\newcommand{\arriba}{5}
\newcommand{\derecha}{8} 			

\newcounter{trunco}

\begin{document}

\title[Minimal $A_\infty$-algebras of endomorphisms]%
{Minimal $A_\infty$-algebras of endomorphisms:\\ The case of $d\ZZ$-cluster
  tilting objects}%

\author[G.~Jasso]{Gustavo Jasso}%

\address[G.~Jasso]{%
  Mathematisches Institut, %
  Universität zu Köln, %
  Weyertal 86-90, %
  50931 Köln, %
  Germany}%
\email{gjasso@math.uni-koeln.de}%
\urladdr{https://gustavo.jasso.info}%

\author[F.~Muro]{Fernando Muro}%

\address[F.~Muro]{%
  Universidad de Sevilla, %
  Facultad de Matemáticas, %
  Departamento de Álgebra, %
  Calle Tarfia s/n, %
  41012 Sevilla, %
  Spain%
}%
\email{fmuro@us.es}%
\urladdr{https://personal.us.es/fmuro/}%



\begin{abstract}
  The Derived Auslander--Iyama Corresponence, a recent result of the authors,
provides a classification up to quasi-isomorphism of the derived endomorphism
algebras of basic $d\ZZ$-cluster tilting objects in $\operatorname{Hom}$-finite
algebraic triangulated categories in terms of a small amount of algebraic data.
In this note we highlight the role of minimal $A_\infty$-algebra structures in
the proof of this result, as well as the crucial role of the enhanced
$A_\infty$-obstruction theory developed by the second-named author.
\end{abstract}

\maketitle

\section*{Introduction}

Let $\kk$ be a ground field. Suppose given a $\operatorname{Hom}$-finite
idempotent-complete algebraic triangulated category $\T$. These assumptions
guarantee that the Krull--Remak--Schmidt Theorem holds in $\T$, see for example~\cite{Kra15}.
Suppose also that $\T$ admits a classical generator $G\in\T$, which for
simplicity we assume to be basic. Hence, $\T$ coincides with its smallest
triangulated subcategory containing $G$ and which is closed under direct
summands. Consider the following naive question:
\begin{quote}
  Can $\T$ be reconstructed---as a triangulated category---from the graded
  endomorphism algebra of the generator $G\in\T$?
\end{quote}
This question is known to have a negative answer. Indeed, let
$\T=\DerCat[b]{\mmod{A}}$ be the bounded derived category of the truncated
polynomial algebra $A=\kk[x]/(x^\ell)$, $\ell\geq 3$. As a consequence of the
Jordan--H{\"o}lder Theorem, $\T$ admits the unique simple $A$-module $S=\kk$ as
a classical generator. A straightforward calculation shows that, independently
from $\ell$, the graded endomorphism algebra of $S$ is
\[
  \Ext[A]{S}{S}\cong\kk[\varepsilon,t]/(\varepsilon^2),\qquad |\varepsilon|=1,\
  |t|=2.
\]
However, the centre $Z(A)=A$ is a derived invariant~\cite{Ric89}, and therefore
$\T$ cannot be reconstructed from the above graded algebra in this case.

A better question is then perhaps the following:
\begin{quote}
  Can $\T$ be reconstructed---as a triangulated category---from the graded
  endomorphism algebra of the generator $G\in\T$ together with some additional
  data?
\end{quote}
Since we assume $\T$ to be algebraic, Keller's Recognition Theorem~\cite{Kel94}
guarantees the existence of a differential graded (=dg) algebra $A$ and an
equivalence of triangulated categories
\[
  \T\stackrel{\sim}{\longrightarrow}\DerCat[c]{A},\qquad G\longmapsto A,
\]
between $\T$ and the perfect derived category of $A$, in particular $\H{A}$ is the graded endomorphism algebra of $G$. Thus, if by `some
additional data' we mean complete knowledge of a derived endomorphism algebra $A$
of $G$ then the latter question has a positive answer. The following issues
arise:
\begin{itemize}
\item Computing the dg algebra $A$ can be highly non-trivial, except in special
  circumstances. For example, if $G$ is a tilting object, then $A$ is the (ungraded) endomorphism algebra of $G$.
\item The quasi-isomorphism class of the dg algebra $A$ is not uniquely
  determined by $\T$, as in principle it depends on a choice of enhancement
  for the latter.
\end{itemize}

We shift the goalposts one last time and ask the following question instead:
\begin{quote}
  Can $\T$ be reconstructed---as a triangulated category---from the graded
  endomorphism algebra of the generator $G\in\T$ together with \emph{a minimal
    amount} of additional data?
\end{quote}
However vague, a general answer the above question seems to be out of reach, at
least without imposing additional assumptions on the generator $G\in\T$. In this
note we report on some of the main results in our joint work~\cite{JKM22}, where
we give an answer to this question when $G$ is a $d\ZZ$-cluster tilting object,
a kind of classical generator that plays an important role in Iyama's higher
Auslander--Reiten Theory~\cite{Iya07,Iya07a,Iya11}.

Following~\cite{IY08,GKO13},
recall that a basic object $G\in\T$ is a \emph{$d\ZZ$-cluster tilting object},
$d\geq1$, if $G\cong\Sigma^d(G)$ and the following are equivalent for an object
$X\in\T$:
\begin{itemize}
\item $X$ is a direct summand of a finite direct sum of copies of $G$.
\item For $0<i<d$, we have $\T(G,\Sigma^i(X))=0$.
\item For $0<i<d$, we have $\T(X,\Sigma^i(G))=0$.
\end{itemize}
Hence, if $d=1$, then $G$ is simply an additive generator of $\T$. Also,
although not immediately obvious from the definition, a $d\ZZ$-cluster tilting
object is necessarily a classical generator of the ambient triangulated category
$\T$. Furthermore, if $G$ is a $d\ZZ$-cluster tilting object, then the graded
endomorphism algebra
\[
  \textstyle\T(G,G)^{\bullet}\coloneqq\bigoplus_{i\in\ZZ}\T(G,\Sigma^i(G))
\]
is concentrated in degrees that are multiples of $d$ and contains an invertible
element of degree $d$, namely a choice of isomorphism $\varphi\colon
G\stackrel{\sim}{\to}\Sigma^d(G)$. Finally, by exploiting the close relationship
between $d\ZZ$-cluster tilting objects and $(d+2)$-angulated categories, it is
shown in \cite{GKO13} that the (ungraded) endomorphism algebra $\T(G,G)$ is a basic
Frobenius algebra. The following is a summary of some of the main results
in~\cite{Mur22,JKM22}; we refer the reader to \cite[Section~6]{JKM22} for
several examples.

\begin{theorem*}[{\cite{Mur22} for $d=1$ and \cite{JKM22} for $d\geq 1$}]\label{thm:1}
  Let $d\geq1$ and assume that the ground field $\kk$ is perfect.\footnote{Or,
    more generally, that the quotient of the finite-dimensional algebra
    $\T(G,G)$ by its Jacobson radical is separable.} Let $\T$ be a
  $\operatorname{Hom}$-finite idempotent-complete algebraic triangulated
  category and ${G\in\T}$ a basic $d\ZZ$-cluster tilting object. Then, up to
  quasi-isomorphism, there is a unique dg algebra $A$ such that there exists
  an equivalence of triangulated categories\footnote{For $d=1$ we need not assume that $\T$ is algebraic, it is so as a consequence.}
  \[
    \T\stackrel{\sim}{\longrightarrow}\DerCat[c]{A},\qquad G\longmapsto A;
  \]
  in particular, $\T$ admits a unique dg enhancement in the sense
  of~\cite{BK90}. Moreover, again up to quasi-isomorphism, the dg algebra $A$
  above can be reconstructed from the basic Frobenius algebra
  $\H[0]{A}\cong\T(G,G)$ and its bimodule
  \[ {\H[-d]{A}\cong\T(G,\Sigma^{-d}(G))}.
  \]
\end{theorem*}

For brevity, in this note we focus on the role of minimal $A_\infty$-structures
in the proof of \Cref{thm:1}, and also in the crucial role played by the
enhanced $A_\infty$-obstruction theory developed by the second-named author
in~\cite{Mur20b} to deal with the case $d=1$. \Cref{thm:1} has a converse, which
associates an algebraic triangulated category with a $d\ZZ$-cluster tilting
object to a twisted $(d+2)$-periodic algebra. Thus, our results provide a
derived version of the Auslader--Iyama Correspondence~\cite{Aus71,Iya07a}, that
can be understood as a classification of the derived endomorphism algebras of
$d\ZZ$-cluster tilting objects in terms of twisted $(d+2)$-periodic algebras. We
refer the reader to the introduction to~\cite{JKM22} for motivation and further
details on this relationship.

\subsection*{An application: The Donovan--Wemyss Conjecture}

The most important application of \Cref{thm:1} to date is, as observed by
Keller in the appendix to~\cite{JKM22}, to the Donovan--Wemyss
Conjecture~\cite{DW16} in the context of the Homological Minimal Model Programme
for threefolds~\cite{Wem18}, see \cite{Wem23} for a recent survey of the rich
algebro-geometric context surrounding the conjecture. \Cref{thm:1} permits to
reconstruct---as a dg-enhanced triangulated category---the singularity category
\[
  \DerCat{R}<sg>\coloneqq\DerCat[b]{\mmod{R}}/\operatorname{K}^{\mathrm{b}}(\operatorname{proj} R)
\]
of an isolated compound Du Val singularity $R$ and assume that it admits a crepant resolution
${p\colon X\to\operatorname{Spec}(R)}$. Associated to $p$ there is the
Donovan--Wemyss contraction algebra
$\Lambda_{\mathrm{con}}=\Lambda_{\mathrm{con}}(p)$, a certain finite-dimensional
algebra that controls the non-commutative deformations of the reduced exceptional fibre of
the resolution. By work of Wemyss, the contraction algebra
$\Lambda_{\mathrm{con}}$ can be realised as the (ordinary) endomorphism algebra
of a $2\ZZ$-cluster tilting object $T$ in the singularity category of $R$, and
the derived endomorphism algebra of $T$ is, according to \Cref{thm:1} and the
$2$-periodicity of the singularity category~\cite{Eis80}, uniquely determined up to 
quasi-isomorphism by the pair $(\Lambda_{\mathrm{con}},\Lambda_{\mathrm{con}})$
consisting of the contraction algebra and its diagonal bimodule. Combined with
results of August~\cite{Aug20a}, Hua--Keller~\cite{HK24} and
Wemyss~\cite{Wem18}, this is enough to finish the proof of the Donovan--Wemyss
Conjecture which, roughly, states that the isomorphism type of the isolated cDV
singularity $R$ is uniquely determined by the derived-equivalence class of the
contraction algebra of any of its crepant resolutions. We refer the reader
to~\cite{JKM22,JKM24} for details on this important application of our main
results. We also mention that an alternative proof of the Dononvan--Wemyss
Conjecture was obtained more recently by Karmazyn, Lepri and
Wemyss~\cite{KLW24}; their proof, which is fundamentally different, also relies
on the use of minimal $A_\infty$-algebra structures and works in an algebraic
setting motivated by the geometric setup in which the conjecture originates.

\section{Hochschild cohomology of graded algebras}

Hochschild cohomology of associative algebras was introduced by the eponymous
author in~\cite{Hoc45}, and its rich algebraic structure was first
investigated by Gerstenhaber in~\cite{Ger63}. In this section we recall the basic
definitions in the case of graded algebras, where sign issues are slightly
more subtle. We adopt the sign conventions used in~\cite{Mur20b}.

Let $\bfLambda$ be a graded algebra.\footnote{For example, the cohomology of the
  dg algebra of endomorphisms of a classical generator of a dg-enhanced
  triangulated category.} The \emph{Hochschild (cochain) complex of $\bfLambda$}
is the bigraded vector space with the components
\[
  \HC[p]<q>{\bfLambda}\coloneqq\Hom[\kk]{\bfLambda^{\otimes p}}{\bfLambda(q)},\qquad
  p\geq 0,\ q\in\ZZ,
\]
where $\bfLambda(q)^i=\bfLambda^{i+q}$, $i\in\ZZ$, and the $\operatorname{Hom}$
is taken in the category of graded vector spaces. It is convenient to introduce
some auxiliary notation in order to describe the differential and the additional
algebraic structure of the Hochschild complex.

Given Hochschild cochains
\[
  c_1\in\HC[p]<q>{\bfLambda}\qquad\text{and}\qquad c_2\in\HC[s]<t>{\bfLambda},
\]
we define, for $1\leq i\leq p$, their \emph{(suspended) infinitesimal composition}
\[
  c_1\bullet_i c_2\in\HC[p+s-1]<q+t>{\bfLambda}
\]
by the formula
\begin{multline*}
  (c_1\bullet_i c_2)(x_1,\dots,x_{p+s-1})\coloneqq\\(-1)^\maltese
  c_1(x_1,\dots,x_{i-1},c_2(x_i,\cdots,x_{i+s-1}),x_{i+s},\cdots x_{p+s-1}),
\end{multline*}
where the sign is given by\footnote{The reader might be surprised that the sign
  does not only involve the Koszul sign. The reason for this is that the use of
  the so-called operadic suspension is implicit,
  see~\cite[Definition~2.4]{Mur16} for the definition with our preferred sign
  conventions.}
\[
  \textstyle\maltese\coloneqq(s-1)(p-i)+t\left(p-1+\sum_{j=1}^{i-1}|x_j|\right).
\]

We also introduce the \emph{pre-Lie product}
\[
  \textstyle\braces{c_1}{c_2}\coloneqq\sum_{i=1}^pc_1\bullet_ic_2
\]
the \emph{Gerstenhaber bracket}
\[
  [c_1,c_2]\coloneqq\braces{c_1}{c_2}-(-1)^{(p+q-1)(s+t-1)}\braces{c_2}{c_1},
\]
and the \emph{cup product}
\[
  \cupp{c_1}{c_2}\coloneqq (\Astr<2>[\bfLambda]\bullet_1
  c_1)\bullet_{p+1}c_2,
\]
where $\Astr<2>[\bfLambda]\in\HC[2]<0>{\bfLambda}$ is the binary multiplication of the algebra $\bfLambda$.
Notice that the pre-Lie product and the Gerstenhaber bracket are bihomogeneous operations of bidegree
$(-1,0)$, and the cup product is bihomogeneous of bidegree $(0,0)$.

The \emph{Hochschild differential} is the bidegree $(1,0)$ map
\[
  \Hd\colon\HC[p]<q>{\bfLambda}\longrightarrow\HC[p+1]<q>{\bfLambda},\qquad x\longmapsto[\Astr<2>[\bfLambda],x].
\]
The cohomology
\[
  \HH{\bfLambda}\coloneqq\H[\bullet,*]{\HC{\bfLambda}}
\]
of the Hochschild complex is the \emph{Hochschild cohomology of $\bfLambda$.}

With respect to the total degree, the Hochschild cohomology of $\bfLambda$ is a
Gerstenhaber algebra, that is it is a graded
commutative algebra which is endowed with a compatible shifted graded Lie algebra
structure. The graded commutative algebra structure is induced
by the cochain-level cup product, while the shifted Lie algebra structure is induced by the
Gerstenhaber bracket.
At the cohomological level, the cup product and the Gerstenhaber bracket satisfy the \emph{Gerstenhaber
  relation},
\[
  [c_1,\cupp{c_2}{c_3}]=\cupp{[c_1,c_2]}{c_3}+(-1)^{(p+q-1)(s+t)}\cupp{c_2}{[c_1,c_3]},
\]
where $c_1\in\HH[p]<q>{\bfLambda}$ and $c_2\in\HH[s]<t>{\bfLambda}$. Moreover, if $\chark{\kk}=2$ then 
$\HH{\bfLambda}$ is equipped with the quadratic operation, called \emph{Gerstenhaber square}, defined at the
cochain level by
\[
  \Sq(c)\coloneqq\braces{c}{c}.
\]
Notice that, if $c$ has
bidegree $(p,q)$, then $\Sq(c)$ has bidegree $(2p-1,2q)$. In general, the
Gerstenhaber square is not $\kk$-linear nor it is an abelian group homomorphism.
Instead, it satisfies the following compatibility relations with respect to the
cup product and the Gerstenhaber bracket:
\begin{align*}
  \label{eq:Gerstehaber_square-relations}
  \Sq(c_1+c_2)&=\Sq(c_1)+\Sq(c_2)+[c_1,c_2],\\
  \Sq(c_1\cdot c_2)&=\Sq(c_1)\cdot c_2^2+c_1\cdot[c_1,c_2]\cdot c_2+c_1^2\cdot\Sq(c_2),\\
  [\Sq(c_1),c_2]&=[c_1,[c_1,c_2]].
\end{align*}
When $\chark{\kk}\neq2$, we can define the Gerstenhaber square in even total degrees by
$\Sq(c)=\tfrac{1}{2}[c,c]$, and the above relations still hold.

\begin{remark}
  The Hochschild cohomology of $\bfLambda$ can be interpreted in terms of the
  bigraded vector space of bimodule self-extensions of the diagonal
  $\bfLambda$-bimodule, that is there are isomorphisms of vector spaces
  \[
    \HH[p]<q>{\bfLambda}\cong\Ext[\bfLambda^e][p,q]{\bfLambda}{\bfLambda}=\Ext[\bfLambda^e][p]{\bfLambda}{\bfLambda(q)},\qquad
    p\geq0,\ q\in\ZZ.
  \]
  Under these identifications, the cup product corresponds to the Yoneda
  product, but the Gerstenhaber bracket has no easy interpretation in these
  terms (see~\cite{Sch98a} for an interpretation of the Gerstenhaber bracket in
  terms of a `loop bracket' of bimodule self-extensions, at least in the case of
  ungraded algebras).
\end{remark}

\begin{remark}
  Together with its algebraic structure, the Hochschild complex is easily seen
  to be functorial on \emph{isomorphisms} of graded algebras.
\end{remark}

\section{Minimal $A_\infty$-algebras and universal Massey products}

Although \Cref{thm:1} is formulated in terms of dg algebras, parts of its proof
rely heavily on the closely-related theory of $A_\infty$-algebras. Standard
references for the latter include~\cite{Kel01,Kel02,Kel02a,Lef03}. Our
presentation makes use of the algebraic structure on the Hochchild cohomology of
a graded algebra, but our sign conventions are precisely those used
in~\cite{Lef03}.

Let $\bfLambda$ be a graded algebra. A \emph{minimal $A_\infty$-algebra
  structure} on $\bfLambda$ is a sequence
\[
  \textstyle\Astr[\bfLambda]=(\Astr<2>[\bfLambda],\Astr<3>[\bfLambda],\dots,\Astr<n>[\bfLambda],\dots)\in\prod_{n\geq
  2}\HC[n]<2-n>{\bfLambda},
\]
where $\Astr<2>[\bfLambda]\in\HC[2]<0>{\bfLambda}$ is the binary product of $\bfLambda$,
that satisfies the \emph{Maurer--Cartan equation}
\[
  0=\braces{\Astr[\bfLambda]}{\Astr[\bfLambda]}.
\]
We call the pair $(\bfLambda,\Astr[\bfLambda])$ a \emph{minimal
  $A_\infty$-algebra}. An \emph{$A_\infty$-morphism}
\[
  f\colon(\bfLambda_1,\Astr[\bfLambda_1])\longrightarrow(\bfLambda_2,\Astr[\bfLambda_2])
\]
between minimal $A_\infty$-algebras is a sequence homogeneous maps
\[
  f_n\colon \bfLambda_1^{\otimes n}\longrightarrow\bfLambda_2,\qquad n\geq 1
\]
of degree $1-n$ that satisfy suitable equations that we do not need to recall
here; it suffices to mention that the \emph{linear part} $f_1$ of $f$ is a
morphism of associative algebras. We call $f$ an \emph{$A_\infty$-isomorphism}
if its linear part is an isomorphism, and a \emph{gauge $A_\infty$-isomorphism}
if it is the identity map (whence necessarily $\bfLambda_1=\bfLambda_2$).

\begin{remark}
  For simplicity, we have omitted the definition of a general
  $A_\infty$-algebra, which is equipped with a degree $1$ differential
  $\Astr<1>[\bfLambda]$ and whose binary operation $\Astr<2>[\bfLambda]$ is only
  associative up to homotopy. There is also a companion notion of
  $A_\infty$-morphism, and one defines the notion
  $A_\infty$-\emph{quasi-}isomorphism in terms of the linear part, which is
  necessarily a morphism of complexes. The advantage of this generality is that
  we may identify dg algebras with $A_\infty$-algebras such that $\Astr<n>[]=0$
  for $n\geq3$, and consider $A_\infty$-morphisms between dg algebras and
  minimal $A_\infty$-algebras.
\end{remark}

\begin{remark}
  Let $(\bfLambda,\Astr[\bfLambda])$ be a minimal $A_\infty$-algebra. The
  Maurer--Cartan equation corresponds to the infinite system of quadratic
  equations
  \[
    \textstyle0=\sum_{p+s-1=n}\braces{\Astr<p>[\bfLambda]}{\Astr<s>[\bfLambda]},\qquad
    n\geq 3,
  \]
  which is one way to write the \emph{$A_\infty$-equations} for an
  $A_\infty$-algebra. For $n=3$, the equation is
  \[
    0=\braces{\Astr<2>[\bfLambda]}{\Astr<2>[\bfLambda]}
  \]
  and corresponds to the associativity of the binary product. More generally,
  the Maurer--Cartan equation expresses the fact that $\Astr[\bfLambda]$
  consists of a family of `coherently associative' higher operations.
\end{remark}

Let $(\bfLambda,\Astr[\bfLambda])$ be a minimal $A_\infty$-algebra. The
$A_\infty$-equation for $n=4$ yields
\[
  0=\braces{\Astr<2>[\bfLambda]}{\Astr<3>[\bfLambda]}+\braces{\Astr<3>[\bfLambda]}{\Astr<2>[\bfLambda]}=[\Astr<2>[\bfLambda],\Astr<3>[\bfLambda]]=\Hd(\Astr<3>[\bfLambda]),
\]
so that $\Astr<3>[\bfLambda]\in\HC[3]<-1>{\bfLambda}$ is a Hochschild cocycle.
Its cohomology class
\[
  \Hclass{\Astr<3>[\bfLambda]}\in\HH[3]<-1>{\bfLambda}
\]
is an invariant of the gauge $A_\infty$-isomorphism class of
$(\bfLambda,\Astr[\bfLambda])$ called the \emph{universal Massey product}.
Moreover, the $A_\infty$-equation for $n=5$ yields
\begin{align*}
  0&=\braces{\Astr<2>[\bfLambda]}{\Astr<4>[\bfLambda]}+\braces{\Astr<4>[\bfLambda]}{\Astr<2>[\bfLambda]}+\braces{\Astr<3>[\bfLambda]}{\Astr<3>[\bfLambda]}\\
   &=[\Astr<2>[\bfLambda],\Astr<4>[\bfLambda]]+\Sq(\Astr<3>[\bfLambda])\\
   &=\Hd(\Astr<4>[\bfLambda])+\Sq(\Astr<3>[\bfLambda]),
\end{align*}
so that
\[
  0=\Sq(\Hclass{\Astr<3>[\bfLambda]})\in\HH[5]<-2>{\bfLambda}.
\]
Finally, set $\Lambda\coloneqq\bfLambda^0$ and let
$j\colon\Lambda\hookrightarrow\bfLambda$ be the canonical inclusion of the
degree $0$ part. The \emph{restricted universal Massey product} is the class
\[
  \textstyle j^*\Hclass{\Astr<3>[\bfLambda]}\in\HH[3]<-1>{\Lambda}[\bfLambda]\cong\Ext[\Lambda^e][3]{\Lambda}{\bfLambda^{-1}},
\]
where the $\operatorname{Ext}$ is computed in the category of ungraded
$\Lambda$-bimodules. The universal Massey product was first investigated
in~\cite{BKS04} (see~\cite{BD89} for its topological precursor, the universal
Toda bracket) and, together with its restricted variant, plays a crucial role in
our work.

\subsection{Minimal models of dg algebras}

Let $A$ be a dg algebra. A \emph{minimal model} for $A$ is a minimal
$A_\infty$-algebra $(\H{A},\Astr[A])$ together with an $A_\infty$-quasi-isomorphism
\[
 i\colon(\H{A},\Astr[A])\stackrel{\sim}{\longrightarrow}A
\]
whose linear part $i_1$ is a \emph{cocycle-selection map}, that is
\[
  \forall x\in\H{A},\quad [i_1(x)]=x.
\]
By a well-known theorem of Kadeishvili~\cite{Kad82}, every dg algebra admits a
minimal model. Since any two minimal models of a dg algebra are related by a
gauge $A_\infty$-isomorphism, we may define the universal Massey product of a dg
algebra $A$ as that of any of its minimal models. The universal Massey product
is an invariant of the quasi-isomorphism class of $A$.

\begin{example}
  \label{ex:counter}
  As in the introduction, let $A=\kk[x]/(x^\ell)$, $\ell\geq 3$, and $S=A/(x)=\kk$. There is
  an isomorphisms of graded algebras
  \[
    \H{\RHom[A]{S}{S}}=\Ext[A]{S}{S}\cong\kk[\varepsilon,t]/(\varepsilon^2),\qquad
    |\varepsilon|=1,\ |t|=2.
  \]
  The dg algebra $\RHom[A]{S}{S}$ admits a minimal model whose only non-zero
  higher operation $\Astr<\ell>[]$ is given by
  \[
    \Astr<\ell>[](\varepsilon,\dots,\varepsilon)=\pm t,
  \]
  see for example~\cite{Mad02,Kel02}.
\end{example}

\begin{remark}
  A dg algebra $A$ is \emph{formal} if it admits a minimal model $(\H{A},\Astr)$
  with vanishing higher operations: $\Astr<n>=0$, $n\geq3$. In general, it is
  difficult to decide whether a given dg algebra is formal but there are useful
  criteria that can actually be checked in practise. For example, a (different)
  theorem of Kadeishvili~\cite{Kad88} (see also~\cite{ST01}) shows that if the Hochschild cohomology of
  $\bfLambda=\H{A}$ vanishes in the range
  \[
    \HH[p+2]<-p>{\bfLambda}=0,\qquad p\geq 1,
  \]
  then every dg algebra with cohomology graded algebra $\bfLambda$ is formal
  (hence $A$ itself is also formal). Similarly, suppose that $B=\kk Q/I$ is a
  bound quiver (finite-dimensional) algebra and let $S$ be the direct sum of a
  complete set of representatives of the simple $A$-modules. Then, the dg
  algebra $A=\RHom[B]{S}{S}$ is formal if and only if $B$ is a Koszul
  algebra~\cite{Kel02}. This also explains why minimal
  $A_\infty$-structures must necessarily appear in \Cref{ex:counter}.
\end{remark}

\begin{remark}
  \label{rmk:quasi-iso-vs-Ai-iso}
  Two dg algebras are quasi-isomorphic if and only if their minimal models are
  $A_\infty$-isomorphic. This means that,
  for our purposes, we may work with minimal $A_\infty$-algebras instead of
  working with dg algebras. The upshot is that minimal $A_\infty$-algebras give
  us access to invariants, such as the universal Massey product, which do not
  seem to be readily available from the dg algebra formalism.
\end{remark}

\subsection{The $d$-sparse case}

Let $d\geq1$ and suppose that the graded algebra $\bfLambda$ is
\emph{$d$-sparse}, that is it is concentrated in degrees that are multiples of
$d$.\footnote{For example, the cohomology of the derived endomorphism algebra of
  a $d\ZZ$-cluster tilting object in a dg-enhanced triangulated category.} Let
$\Astr[\bfLambda]$ be a minimal $A_\infty$-algebra structure on $\bfLambda$. By
degree reasons, we must have $\Astr<n>[\bfLambda]=0$ if $2-n\not\in d\ZZ$. Thus,
the only possibly non-zero higher operations are
\[
  \Astr<d+2>[\bfLambda],\ \Astr<2d+2>[\bfLambda],\ \Astr<3d+2>[\bfLambda],\
  \dots,\ \Astr<kd+2>[\bfLambda],\ \dots
\]
The $A_\infty$-equation for $n=d+3$ yields the fact that
$\Astr<d+2>[\bfLambda]\in\HC[d+2]<-d>{\bfLambda}$ is a Hochschild cocycle, whose
class
\[
  \Hclass{\Astr<d+2>[\bfLambda]}\in\HH[d+2]<-d>{\bfLambda}
\]
we call the \emph{universal Massey product of length $d+2$}. The
$A_\infty$-equation for $n=2d+3$ yields the vanishing of its Gerstenhaber
square:
\[
  0=\Sq(\Hclass{\Astr<d+2>[\bfLambda]})\in\HH[2d+3]<-2d>{\bfLambda}.
\]
Due to the $d$-sparseness condition, this variant of the universal Massey product
is also an invariant of the gauge $A_\infty$-isomorphism class of
$(\bfLambda,\Astr[\bfLambda])$. Let $j\colon\Lambda\hookrightarrow\bfLambda$ be
the canonical inclusion of the degree $0$ part. The \emph{restricted universal
  Massey product of length $d+2$} is the class
\[
  \textstyle j^*\Hclass{\Astr<d+2>[\bfLambda]}\in\HH[d+2]<-d>{\Lambda}[\bfLambda]\cong\Ext[\Lambda^e][d+2]{\Lambda}{\bfLambda^{-d}},
\]
where the $\operatorname{Ext}$ is computed in the category of ungraded
$\Lambda$-bimodules.

If $A$ is a dg algebra with $d$-sparse cohomology, we define its universal
Massey product of length $d+2$ as that of any of its minimal models, and this
yields an invariant of the quasi-isomorphism class of $A$. 

\section{The Derived Auslander--Iyama Correspondence}

The following is one of the main results in~\cite{JKM22}; it implies \Cref{thm:1}.

\begin{theorem}[{\cite[part of Theorem~A and Corollary~4.5.17]{JKM22}}]
  Suppose that the ground field $\kk$ is perfect.
  \label{thm:A-injectivity}
  Let $d\geq 1$ and $A$ a dg algebra whose cohomology is $d$-sparse; set
  \[
    \bfLambda\coloneqq\H{A}.
  \]
  Suppose that $\Lambda\coloneqq\bfLambda^0$ is a basic Frobenius algebra and
  that there exists an invertible element of degree $d$ in $\bfLambda$. The following
  statements are equivalent:
  \begin{enumerate}
  \item The free dg $A$-module $A\in\DerCat[c]{A}$ is a basic $d\ZZ$-cluster
    tilting object.
  \item The restricted universal Massey product of length $d+2$ of $A$,
    \[
      \textstyle j^*\Hclass{\Astr<d+2>[A]}\in\HH[d+2]<-d>{\Lambda}[\bfLambda]\cong\Ext[\Lambda^e][d+2]{\Lambda}{\bfLambda^{-d}},
    \]
    represents an isomorphism
    $\bfLambda^{-d}\cong\Omega_{\Lambda^e}^{d+2}(\Lambda)$ in the stable
    category of $\Lambda$-bimodules.
  \end{enumerate}
  Moreover, if the above conditions are satisfied, then $A$ is uniquely
  determined up to quasi-isomorphism by the pair $(\Lambda,\bfLambda^{-d})$
  consisting of the basic Frobenius algebra $\Lambda$ and the
  (invertible\footnote{Indeed, the multiplication maps
    \[\bfLambda^{d}\otimes_{\Lambda}\bfLambda^{-d}\stackrel{\sim}{\longrightarrow}\Lambda\qquad\text{and}\qquad\bfLambda^{-d}\otimes_{\Lambda}\bfLambda^{d}\stackrel{\sim}{\longrightarrow}\Lambda\] 
  are invertible by assumption.
  })
  $\Lambda$-bimodule $\bfLambda^{-d}$.
\end{theorem}

\begin{remark}
  Part of the proof of \Cref{thm:A-injectivity} relies on a careful analysis of the
  relationship between restricted universal Massey products, the standard
  $(d+2)$-angulated category structure on the additive closure of a
  $d\ZZ$-cluster tilting object in a triangulated category~\cite{GKO13}, and
  Amiot--Lin $(d+2)$-angulations~\cite{Ami07,Lin19}.
\end{remark}

\Cref{thm:A-injectivity} has the following converse.

\begin{theorem}[{\cite[Theorem~5.1.2]{JKM22}}]
  \label{thm:A-surjectivity}
  Let $d\geq 1$. Let $\Lambda$ be a basic Frobenius algebra and
  $\sigma\colon\Lambda\stackrel{\sim}{\to}\Lambda$ an algebra automorphism. Form the
  $\sigma$-twisted Laurent polynomial algebra\footnote{Up to isomorphism of
    graded algebras, $\bfLambda(\sigma,d)$ only depends on the class of $\sigma$
  in the outer automorphism group of $\Lambda$.}
  \[
    \bfLambda=\bfLambda(\sigma,d)\coloneqq\frac{\Lambda\langle\imath^\pm\rangle}{(\imath
      x-\sigma(x)\imath)_{x\in\Lambda}},\qquad |\imath|=-d.
  \]
  Suppose that there exists a Yoneda class
  \[
    \textstyle\eta\in\Ext[\Lambda^e][d+2]{\Lambda}{\bfLambda^{-d}}\cong\HH[d+2]<-d>{\Lambda}[\bfLambda]
  \]
  that represents an isomorphism
  $\bfLambda^{-d}\cong\Omega_{\Lambda^e}^{d+2}(\Lambda)$ in the stable category
  of $\Lambda$-bimodules. Then, up to gauge $A_\infty$-isomorphism, there exists
  a unique minimal $A_\infty$-algebra structure $\Astr[\eta]$ on $\bfLambda$
  whose restricted universal Massey product of length $d+2$ satisfies
  \[
    j^*\Hclass{\Astr<d+2>[\eta]}=\eta\in\Ext[\Lambda^e][d+2]{\Lambda}{\bfLambda^{-d}}.
  \]
\end{theorem}

\begin{remark}
  In \cite{JKM22}, we have formulated
  \Cref{thm:A-injectivity,thm:A-surjectivity} as a single result, namely as
  Theorem~A in~\emph{loc.~cit.}, which better illustrates the analogy with the
  Auslander--Iyama Correspondence~\cite{Aus71,Iya07a}. In these proceedings we
  have opted for the above alternative formulation of our main results as it
  highlights the role of universal Massey products in our work.
\end{remark}

\section{Enhanced $A_\infty$-obstruction theory}

The proofs of \Cref{thm:A-injectivity,thm:A-surjectivity} rely crucially on the
enhanced $A_\infty$-obstruction theory developed by the second-named author
in~\cite{Mur20b}. In this section we explain how this obstruction theory arises,
and why it is useful to leverage it in this context. For simplicity we only
consider the case $d=1$ since, surprisingly enough, the general $d$-sparse case
is not much more complicated.

Obstruction theory has its origin in algebraic topology where it was first used
to solve problems such as the following. Let $W$ and $X$ be CW complexes with $X$ simply connected and $A\subset
W$ a subcomplex. Suppose given a continuous map $f\colon A\to
X$. One wishes to answer
the following question:
\begin{quote}
  When does there exist an extension $g\colon W\to X$ of $f$, that is a morphism
  rendering the diagram below commutative?
  \[
    \begin{tikzcd}
      A\dar[hookrightarrow]\rar{f}&X\\
      W\urar[dotted,swap]{g}
    \end{tikzcd}
  \]
\end{quote}
By standard results, one can describe $X\simeq\varprojlim X_n$ as the (homotopy)
limit of a Postnikov tower of principal fibrations (the precise definitions are
not relevant for the discussion)
\[
  \cdots\twoheadrightarrow X_n\twoheadrightarrow\cdots \twoheadrightarrow
  X_2\twoheadrightarrow X_1\twoheadrightarrow X_0=*
\]
The original extension problem can then be solved by induction going up along
the tower. Suppose given a commutative square
\[
  \begin{tikzcd}
    A\dar[hook]\rar{f}\ar[phantom]{dr}[description]{=}&X\dar[two heads]\\
    W\rar[swap]{g_{n-1}}&X_{n-1}
  \end{tikzcd}
\]
In the inductive step one needs to find a lift
\[
  \begin{tikzcd}
    A\dar[hook]\rar{f}&X\rar&X_n\dar\\
    W\ar[dotted]{urr}\ar{rr}[swap]{g_{n-1}}&&X_{n-1}
  \end{tikzcd}
\]
that is a map $W\to X_n$ rendering the diagram commutative. To the latter
lifting problem one associates an \emph{obstruction class}
\[
  \omega_n\in\H[n+1]{W,A;\pi_n(X)}
\]
in the relative cohomology of the pair $(W,A)$ with coefficients in the $n$-th
homotopy group $\pi_n(X)$. The punchline is that the required lift $W\to X_n$
exists if and only if the obstruction class vanishes: $\omega_n=0$. In
particular, if for all $n\geq1$ one has
\[
  {\H[n+1]{W,A;\pi_n(X)}=0},
\]
then the original extension problem can always be solved since all possible obstructions
necessarily vanish, see~\cite[Sec.~4.3]{Hat02} for a textbook treatment of this
problem.

Thus, given a suitable existence problem, an obstruction theory associates a
class in a suitable cohomology theory whose vanishing corresponds precisely to
the existence of a solutions to the problem. The situation is also familiar in
algebra: Let $A$ be an algebra and $M$ and $N$ a pair of $A$-modules. Suppose given a submodule
$L\subset M$ and a morphism of $A$-modules $f\colon L\to N$. Then, the
obstruction to the existence of an extension of $f$ to $M$,
\[
  \begin{tikzcd}
    L\rar{f}\dar[hookrightarrow]&N\\
    M\urar[dotted]
  \end{tikzcd}
\]
is precisely the class
\[ [f_*\delta]\in\Ext[A][1]{M/L}{N}
\]
determined by the push forward along $f$ of the induced short exact sequence
\[
  \delta\colon0\to L\hookrightarrow M\twoheadrightarrow M/L\to0.
\]

We wish to recast \Cref{thm:A-surjectivity} in terms of obstruction theory. For
this, it is convenient to interpret minimal $A_\infty$-structures on a graded
algebra $\bfLambda$ as morphisms in a suitable category, namely in the category
of dg operads. Roughly speaking, a \emph{dg operad} $\O$ is a sequence $\O(n)$,
$n\geq0$, of dg vector spaces (=cochain complexes of vector spaces) equipped
with \emph{infinitesimal compositions}
\[
  \circ_i\colon\O(p)\otimes\O(q)\longrightarrow\O(p+q-1),\qquad 1\leq i\leq
  p,\quad q\geq0,
\]
that are associative and unital in a suitable sense. The prototypical example of
a dg operad is the \emph{endomorphism dg operad} $\opEnd{V}$ of a dg vector
space $V$, whose dg vector spaces of operations are
\[
  \opEnd{V}\!(n)\coloneqq\dgHom[\kk]{V^{\otimes n}}{V},\qquad n\geq0,
\]
and with the infinitesimal compositions given by the apparent compositions of
multilinear maps. There is also a dg operad $\opA$ with the property that
morphisms of dg operads $\opA\to\opEnd{V}$ are in canonical bijection with
$A_\infty$-structures on $V$. Moreover, the dg operad $\opA$ admits an
exhaustive filtration by dg suboperads
\[
  \opA[2]\hookrightarrow\opA[3]\hookrightarrow\cdots\opA[n]\hookrightarrow\cdots\hookrightarrow\opA,
\]
where the dg operad $\opA[n]$ is such that morphisms of dg operads
$\opA[n]\to\opEnd{V}$ are in canonical bijection with \text{$A_n$-structures} on
$V$, that is truncated $A_\infty$-stuctures with operations only up to arity
$n$.

The category of dg operads is equipped with a Quillen model category structure,
namely the projective model structure transferred from the category of dg vector
spaces, see~\cite{Lyu11,Mur11}. Consequently, there is a notion of homotopy
between morphisms of dg operads (all dg operads are fibrant and the dg operads
$\opA[n]$, ${2\leq n\leq \infty}$, are examples of cofibrant objects in this
model structure), and it turns out that two morphisms of dg operads
$\opA[n]\to\opEnd{V}$ are homotopic if and only if their corresponding
$A_n$-structures are gauge $A_n$-isomorphic. Therefore we can use techniques of
abstract homotopy theory for investigating $A_\infty$-algebra structures.

Let $\bfLambda$ be a graded algebra. Suppose given a morphism of dg operads
${f\colon\opA[4]\to\opEnd{\bfLambda}}$ corresponding to a minimal
$A_4$-structure $(\Astr<2>[f],\Astr<3>[f],\Astr<4>[f])$ on $\bfLambda$ with
$\Astr<2>[f]=\Astr<2>[\bfLambda]$. Consider the following question:
\begin{quote}
  When does there exist a morphism of dg operads
  $g\colon\opA[5]\to\opEnd{\bfLambda}$ such that the restriction of $f$ and $g$
  to $\opA[3]$ coincide?
  \[
    \begin{tikzcd}
      \opA[3]\dar[hookrightarrow]\rar[hookrightarrow]\ar[phantom]{dr}[description]{=}&\opA[4]\dar{f}\\
      \opA[5]\rar[dotted,swap]{g}&\opEnd{\bfLambda}
    \end{tikzcd}
  \]
\end{quote}
Firstly, the fact that $f$ determines a minimal $A_4$-structure on $\bfLambda$
implies that the corresponding operation $\Astr<3>[f]\in\HC[3]<-1>{\bfLambda}$
is a Hochchild cocycle (we have seen that this is a consequence of the
$A_\infty$-equation for $n=4$). It follows from the $A_\infty$-equation for
$n=5$ that the obstruction class corresponding to this extension problem is
precisely the Gerstenhaber square
\[
  \Sq(\Hclass{\Astr<3>[f]})\in\HH[5]<-2>{\bfLambda}.
\]
Hence, this class vanishes if and only if there exists a morphism
${g\colon\opA[5]\to\opEnd{\bfLambda}}$ with the required property. Notice that
such a morphism $g$ determines a minimal $A_5$-structure
\[
  (\Astr<2>[g],\Astr<3>[g],\Astr<4>[g],\Astr<5>[g])
\]
on $\bfLambda$ such that
\[
  (\Astr<2>[f],\Astr<3>[f])=(\Astr<2>[g],\Astr<3>[g]).
\]
In other words, precisely when $\Sq{\Hclass{\Astr<3>[f]}}=0$, we may extend the
given minimal $A_4$-structure $\Astr[f]$ on $\bfLambda$ to a minimal
$A_5$-structure $\Astr[g]$, possibly after replacing the operation $\Astr<4>[f]$
but, and this is crucial, without replacing $\Astr<3>[f]$ and hence fixing the
universal Massey product
\[
  \Hclass{\Astr<3>[f]}\in\HH[3]<-1>{\bfLambda}.
\]
In principle, it is possible to continue analysing such extension problems in
order to determine whether it is possible to extend the minimal $A_4$-structure
$\Astr[f]$ to a full minimal $A_\infty$-structure on $\bfLambda$ with identical
universal Massey product (compare with the existence statement
in~\Cref{thm:A-surjectivity}). However, it is more convenient to first
understand a little bit better the nature of our problem.

By general results of Dwyer and Kan~\cite{DK80b} for Quillen model categories,
given two dg operads $\O$ and $\P$, with $\O$ cofibrant, there is a space (=Kan complex)
\[
  \Map{\O}{\P}
\]
whose vertices are the 
morphisms of dg operads $\O\to\P$. It is not difficult to see that 
$\Map{\opA}{\opEnd{\bfLambda}}$ is the (homotopy) limit of the induced tower (of Kan
fibrations)
\[
  \cdots\twoheadrightarrow\Map{\opA[n]}{\opEnd{\bfLambda}}\twoheadrightarrow \cdots\twoheadrightarrow\Map{\opA[3]}{\opEnd{\bfLambda}}\twoheadrightarrow\Map{\opA[2]}{\opEnd{\bfLambda}}.
\]
For brevity, let us write
\[
  X_n\coloneqq\Map{\opA[n+2]}{\opEnd{\bfLambda}},\qquad 0\leq n\leq \infty.
\]
One may hope to understand the space $X_\infty$ in terms of
the above tower. For example, given a minimal
$A_\infty$-structure $\Astr[f]$ on $\bfLambda$ defined by a point $f\in X_\infty$, for each $i\geq 0$ there is a
Milnor short exact sequence~\cite[Section~IX.3.1]{BK72}
\[
  *\to\varprojlim\nolimits^1
  \pi_{i+1}(X_n)\to\pi_i(X_\infty)\to\varprojlim\pi_i(X_n)\to *,
\]
of pointed sets for $i=0$ and of groups for $i>0$,
where $\varprojlim^1$ is the first derived functor of the limit functor
(see~\cite[Section~IX.2]{BK72} for the definition in the case of arbitrary
groups) and all homotopy groups are based at the points induced by $f$. For
$i=0$, the above sequence can be used to determine when two minimal
$A_\infty$-structures on $\bfLambda$ are gauge $A_\infty$-isomorphic (compare
with the uniqueness statement in~\Cref{thm:A-surjectivity}). For this, it would
be useful to have a criterion for the vanishing of the term
$\varprojlim\nolimits^1 \pi_{i+1}(X_n)$. Explaining this criterion requires
more machinery.

We fix a minimal
$A_\infty$-structure on $\bfLambda$ given by a point $f\in X_\infty$. The restrictions of $f$ determine minimal
$A_n$-structures so that we can regard
\[
  \cdots\twoheadrightarrow X_n\twoheadrightarrow\cdots\twoheadrightarrow
  X_2\twoheadrightarrow X_1\twoheadrightarrow X_0
\]
as a pointed tower of spaces. To such a tower, Bousfield and Kan~\cite{BK72}
associate a spectral sequence $\BK[r][s]<t>$. This spectral sequence is somewhat
anomalous: it is defined only in the region $0\leq s\leq t$ and, moreover, some
of its terms are only groups and others are only pointed sets (rather than
abelian groups). Since spectral sequences might be unfamiliar to many of our
target readers, we shall spend some time explaining why it is useful to consider
such a spectral sequence in this context (a standard reference for the theory of
spectral sequences and its applications is~\cite{McC01}). In fact, as we explain below, in order
to place the obstructions needed to solve the kind of extension problems posed
by \Cref{thm:A-surjectivity}, an extension (no pun intended) of the
Bousfield--Kan spectral sequence constructed by the second author
in~\cite{Mur20b} is essential.

We begin by reviewing the most relevant aspects of the construction of the
aforementioned Bousfield--Kan spectral sequence. For $r\geq1$ and $0\leq s\leq
t$, the term $\BK[r][s]<t>$ is given by the homology of the sequence
\[
  \resizebox{0.975\hsize}{!}{
    $\operatorname{Ker}[\pi_{t-s+1}(X_{s-1})\to\pi_{t-s+1}(X_{s-r})]\to\pi_{t-s}(X_{s,s-1})\to\frac{\pi_{t-s}(X_s)}{\operatorname{Im}[\pi_{t-s}(X_{s+r-1})\to\pi_{t-s}(X_s)]},$}
\]
where $X_{s,s-1}$ is the (homotopy) fibre of the map $X_s\twoheadrightarrow
X_{s-1}$ (when $s=t$ one must rather take the quotient of a natural action of
the group on the left on the pointed kernel of the map on the right). The reader
need not focus on the precise definition of the term $\BK[r][s]<t>$, but only
notice that the definition depends on the homotopy groups of the terms of the
tower and of its consecutive fibres. The Bousfield--Kan spectral sequence has
the following properties, which reflect the structure present on the homotopy
groups of a pointed space:
\begin{itemize}
\item For all $0\leq s\leq t$, the term $\BK[r][s]<t>$ is a pointed set.
\item The pointed set $\BK[r][s]<t>$ is a group if $t-s\geq 1$.
\item The group $\BK[r][s]<t>$ is abelian for $t-s\geq 2$.
\end{itemize}
For fixed $r\geq1$, the terms $\BK[r]$ form what is called the \emph{$r$-th
  page} of the spectral sequence. In addition, for fixed $r\geq1$ there are
bidegree $(r,r-1)$ differentials
\[
  d_r\colon\BK[r][s]<t>\longrightarrow\BK[r][s+r]<t+r-1>.
\]
Of course, the differentials are only defined when both the source and the
target are defined (their precise definition is not needed in our discussion).
We also remind the reader of the key defining property of a spectral sequence,
namely that the terms of $(r+1)$-page are given by the cohomology of the $r$-th
page whenever the formula makes sense:\footnote{In this formula we adopt the
  convention that for $s<r$ the incoming differential is $0\to\BK[r][s]<t>$.}
\[
  \H[s,t]{\BK[r]}\cong\BK[r+1][s]<t>.
\]
Notice that the differential of the $(r+1)$-page is defined independently of the
previous pages of the spectral sequence.

By the Complete Convergence Lemma in \cite[Section~IX.5.4]{BK72}, if for fixed
$i\geq 0$ and each $s\geq0$ we have
\[
  \varprojlim\nolimits^1\BK[r][s]<s+i>=*,
\]
then also
\[
  \varprojlim\nolimits^1\pi_{i+1}(X_n)=*.
\]
This is the criterion that we announced a couple of paragraphs above and it can
be verified in practise. However, the Bousfield--Kan spectral sequence does not
contain enough information to apply this criterion in combination with the
Milnor short exact sequence to deduce the kind of statements needed to
prove~\Cref{thm:A-surjectivity}.

It is clear from the construction why the Bousfield--Kan spectral sequence is
only defined in the region $0\leq s\leq t$: a space only has non-negative
homotopy groups. This elementary observation also suggests that, in order to
extend the region on which the Bousfiled--Kan spectral sequence is defined, one
could try to realise at least some of the spaces involved in the
definition\footnote{The definition of the differentials in the Bousfield--Kan
  spectral sequence makes use of the homotopy fibres of the maps
  $X_n\twoheadrightarrow X_m$, $n\geq m$.} as underlying spaces of suitable
spectra in the sense of stable homotopy theory (see~\cite{Gre07} for a very
readable survey on spectra). Indeed, if $X\simeq\Omega^{\infty}(\widetilde{X})$
is the underlying space of a spectrum $\widetilde{X}$, then there are
isomorphisms
\[
  \pi_n(\widetilde{X})\cong\pi_n(X),\qquad n\geq0,
\]
but the spectrum $\widetilde{X}$ can have non-zero negative homotopy groups!

\setcounter{rpage}{5}
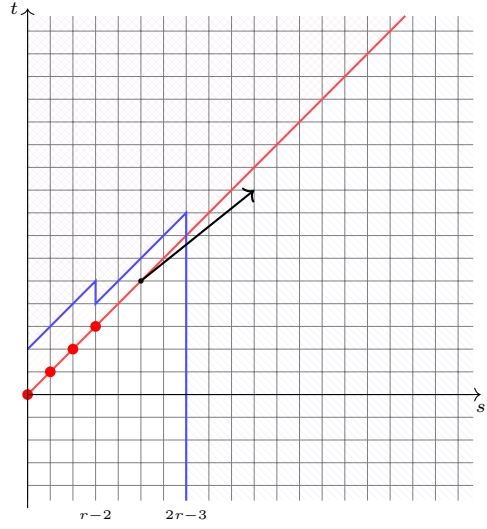
\begin{figure}[ht]
	\begin{tikzpicture}
	\draw[step=\rescale,gray,very thin] (0,{-\abajo*\rescale+\margen}) grid ({(3*\rpage -3 +\derecha)*\rescale-\margen},{(3*\rpage -3 +\arriba)*\rescale-\margen});
	\filldraw[pattern=north east lines, pattern color=red,,draw=none,opacity=0.2] (0,0) -- ({(3*\rpage -3 +\arriba)*\rescale-\margen},{(3*\rpage -3 +\arriba)*\rescale-\margen}) --
	(0,{(3*\rpage -3 +\arriba)*\rescale-\margen}) -- cycle;
	\draw[red!70,thick]  (0,0) -- ({(3*\rpage -3 +\arriba)*\rescale-\margen},{(3*\rpage -3 +\arriba)*\rescale-\margen});
	\foreach  \x in {2,...,\rpage}
	\node[fill=red,draw=none,circle,inner sep=.5mm,opacity=1]   at ({(\x-2)*\rescale},{(\x-2)*\rescale}) {};
	\filldraw[pattern=north west lines, pattern color=blue,draw=none,opacity=0.2]  (0,{(3*\rpage -3 +\arriba)*\rescale-\margen}) --
	(0,2*\rescale) -- ({(\rpage-2)*\rescale},{\rpage*\rescale}) -- ({(\rpage-2)*\rescale},{(\rpage-1)*\rescale}) --({(2*\rpage-3)*\rescale},{(2*\rpage-2)*\rescale}) -- ({(2*\rpage-3)*\rescale},{-\abajo*\rescale+\margen}) -- ({(3*\rpage -3 +\derecha)*\rescale-\margen},{-\abajo*\rescale+\margen}) --
	({(3*\rpage -3 +\derecha)*\rescale-\margen},{(3*\rpage -3 +\arriba)*\rescale-\margen}) -- cycle;
	\draw[blue!70,thick] (0,2*\rescale) -- ({(\rpage-2)*\rescale},{\rpage*\rescale}) -- ({(\rpage-2)*\rescale},{(\rpage-1)*\rescale}) --({(2*\rpage-3)*\rescale},{(2*\rpage-2)*\rescale}) -- ({(2*\rpage-3)*\rescale},{-\abajo*\rescale+\margen}) node[black,anchor=north] {$\scriptscriptstyle 2r-3$};
	\node[black,anchor=north] at ({(\rpage-2)*\rescale},{-\abajo*\rescale+\margen}) {$\scriptscriptstyle r-2$};
	\draw [->] (0,0)  -- ({(3*\rpage -3 + \derecha)*\rescale},0) node[anchor=north] {$\scriptstyle s$};
	\draw [->] (0,-{\abajo*\rescale}) -- (0,{(3*\rpage -3 +\arriba)*\rescale})node[anchor=east] {$\scriptstyle t$};
	\node[fill=,draw=none,circle,inner sep=.25mm,opacity=1]   at ({\rpage*\rescale},{\rpage*\rescale}) {};
	\draw [black,thick,->] ({\rpage*\rescale},{\rpage*\rescale})  -- ({2*\rpage*\rescale},{(2*\rpage-1)*\rescale});
	\end{tikzpicture}
	\caption{Range of definition of the extended Bousfield--Kan spectral sequence
    with $r=\arabic{rpage}$. The region hashed NE-SW consists of vector spaces, the
    region hashed SE-NW consists of abelian groups, with the exception of the large dots that
    are plain pointed sets. We depict a differential that `jumps' from the range
    of definition of the classical Bousfield--Kan spectral sequence (hashed SE-NW)
    to its extended part.}
  \label{fig:mrd}
\end{figure}

Following the above heuristics, in~\cite{Mur20b} the second-named author
constructed a specific extension of the Bousfield--Kan spectral sequence of the
tower of spaces of minimal $\opA[n]$-structures on a graded vector space. The first task
that one typically performs with a spectral sequence is that of computing its
first page:
\begin{itemize}
\item The terms $\BK[r][s]<t>$ are defined in the right half-plane, except if
  $t<s<2r-3$, see \Cref{fig:mrd}.
\item The first page of the extended Bousfield--Kan spectral is given by the
  Hochshcild cochains of $\bfLambda$:
  \[
    \BK[1][s]<t>=\HC[s+2]<-t>{\bfLambda},\qquad s\geq0,\quad t\in\ZZ.
  \]
\item The differentials on the first page are given by the Hochschild
  differentials:
  \begin{align*}
    d_1\colon\BK[1][s]<t>&\longrightarrow\BK[1][s+1]<t>,\qquad
                           s\geq 0,\quad t\in\ZZ,\\
    c&\longmapsto[\Astr<2>[f],c].
  \end{align*}
  Thus, the first page of the extended Bousfield--Kan spectral sequence can be
  identified with part of the Hochschild complex $\HC{\bfLambda}$.
\item As a consequence of the previous two bullet points, the second page is
  (mostly) given by the Hochschild cohomology of $\bfLambda$:
  \[
    \BK[2][s]<t>=\HH[s+2]<-t>{\bfLambda},\qquad s>0,\quad t\in\ZZ.
  \]
\end{itemize}
After computing the first page, and with that the terms of the second page, the
next task is that of computing the differential of the second page, and in
general this is non trivial. In this case the result is the following:
\begin{itemize}
\item The differentials on the second page are (mostly) given by the
  Gerstenhaber bracket with the universal Massey product:
  \begin{align*}
    d_2\colon\BK[2][s]<t>&\longrightarrow\BK[2][s+2]<t+1>,\qquad
                           s\geq 2,\quad t\in\ZZ,\\
    c&\longmapsto[\Hclass{\Astr<3>[f]},c].
  \end{align*}
\end{itemize}
This is an important computation that motivates the introduction of the
Hochschild--Massey complex of the pair $(\bfLambda,\Hclass{\Astr<3>[f]})$
in~\cite[Definition~5.2.5]{JKM22}, see also Theorem~B therein.

Next, we describe some other terms of the extended Bousfield--Kan
spectral sequence in a way that its usefulness in our context becomes more
apparent.
\begin{itemize}
\item For $0\leq s\leq r-1$, the term $\BK[r][s]<s>$ is the pointed set of
  homotopy classes of morphisms of dg operads
  $g\colon\opA[s+2]\to\opEnd{\bfLambda}$ that can be extended to $\opA[s+r+1]$
  and whose restriction to $\opA[s+1]$ is homotopic to the restriction of $f$:
  \[
    \begin{tikzcd}
      \opA[s+1]\dar\rar\ar[phantom]{dr}[description]{\Rightarrow}&\opA[s+2]\dar{g}\rar&\opA[s+r+1]\dlar[dotted]\\
      \opA\rar[swap]{f}&\O
    \end{tikzcd}
  \]
  The base point is the restriction of $f$ to $\opA[s+2]$. For example, if
  $\BK[r][s]<s>=*$, then every morphism of dg operads $g$ as above has the same restriction to $\opA[s+1]$ as $f$, up to homotopy. This can be leveraged when
  considering the Milnor short exact sequence of pointed sets
  \[
    *\to\varprojlim\nolimits^1\pi_{1}(X_n)\to\pi_0(X_\infty)\to\varprojlim\pi_0(X_n)\to
    *,
  \]
  where we remind the reader that $\pi_0(X_n)$ is the set of homotopy classes of
  maps $\opA[n]\to\opEnd{\bfLambda}$ or, equivalently, the set of minimal
  $A_n$-structures on $\bfLambda$ up to gauge $A_n$-isomorphism.
\item There is a pointed bijection between the pointed set of homotopy classes
  of maps $g\colon\opA[3]\to\opEnd{\bfLambda}$ which extend to $\opA[4]$ and
  such that the restrictions of $f$ and $g$ to $\opA[2]$ are homotopic (this last condition says that $\Astr<2>[g]= \Astr<2>[f] = \Astr<2>[\bfLambda] $ is the binary product of the graded algebra $\bfLambda$),
  \[
    \begin{tikzcd}
      \opA[2]\dar\rar\ar[phantom]{dr}[description]{\Rightarrow}&\opA[3]\dar{g}\rar&\opA[4]\dlar[dotted]\\
      \opA[\infty]\rar[swap]{f}&\opEnd{\bfLambda},
    \end{tikzcd}
  \]
  and the term $\BK[2][1]<1>=\HH[3]<-1>{\bfLambda}$. The base point in the
  source is the restriction of $f$ to $\opA[3]$ and the base point in the target
  is $0$. The bijection maps the homotopy class of
  $g\colon\opA[3]\to\opEnd{\bfLambda}$ to
  $\Hclass{\Astr<3>[g]}-\Hclass{\Astr<3>[f]}$. Consequently, the maps $f$ and
  $g$ induce identical universal Massey products if and only if their
  restrictions to $\opA[3]$ are gauge $A_3$-isomorphic.
\end{itemize}
Finally, the following observation is also fundamental:
\begin{itemize}
\item If we fix a minimal $A_n$-structure on $\bfLambda$ given by a map $f\colon\opA[n]\to\opEnd{\bfLambda}$ for some $n<\infty$, then
  the extended Bousfield--Kan spectral sequence can be defined up to a certain
  page (determined by an increasing function of $n$). This is important when
  studying extensions of minimal $A_n$-structures on $\bfLambda$ to a full
  minimal $A_\infty$-structure (compare with
  \Cref{thm:A-surjectivity}).
\end{itemize}

The above properties of the extended Bousfield--Kan spectral sequence are crucial to
the proof of~\cite[Theorem~B]{JKM22}, which is itself a key ingredient in the proof
of \Cref{thm:A-surjectivity}. In a nutshell, this spectral sequence provides an
approach to studying the vanishing of obstructions to the existence and
uniqueness of minimal $A_\infty$-structures on a graded algebra.

\begin{remark}
  Another question that is usually dealt with when studying a spectral
  sequence is that of convergence. Since in the limit $r\to\infty$ of the
  Bousfield--Kan spectral sequence and its extension agree, one knows that the
  term $\BK[\infty][s]<t>$, ${t-s>0}$, contributes to the homotopy group
  $\pi_{t-s}(X_\infty)$ of the homotopy limit of the tower and to the set of
  connected components $\pi_0(X_\infty)$ when $s=t$,
  see~\cite[Section~IX.5]{BK72}. In the case of the space
  \[
    X_\infty=\Map{\opA}{\opEnd{\bfLambda}}
  \]
  of minimal $A_\infty$-algebra structures on the underlying graded vector
  space of $\bfLambda$, these are given as follows~\cite[Proposition~6.9]{Mur20b}:
  \begin{align*}
    \pi_0(X_\infty)&=\text{minimal $A_\infty$-algebras with underlying graded vector space $\bfLambda$}\\ 
    &\phantom{=}\text{\ \ modulo gauge $A_\infty$-isomorphisms,}\\
    \pi_1(X_\infty,\Astr[\bfLambda])&=\text{gauge $A_\infty$-automorphisms of the minimal $A_\infty$-algebra $(\bfLambda,\Astr[\bfLambda])$}\\&\phantom{=}\text{\ \ modulo homotopies with identity linear part,}\\
    \pi_n(X_\infty,\Astr[\bfLambda])&=\operatorname{HH}^{2-n}_{\geq2}(B)\qquad n\geq2,
  \end{align*}
  where we write $B$ for the minimal $A_\infty$-algebra $(\bfLambda,\Astr[\bfLambda])$ and
  $\operatorname{HH}^{\bullet}(B)$ for its Hochschild cohomology,
  see \cite{Mur20b} for definitions and explanations of the notation.
\end{remark}

\begin{remark}
  Based heavily on the above ideas, in \cite{JM25} we develop an obstruction
  theory for the existence and uniqueness of minimal \emph{$A_\infty$-bimodule}
  structures. The main results in \cite{JM25} are then used in \cite{JM25b} to
  prove a `Calabi--Yau variant' of the Derived Auslander--Iyama Correspondence.
\end{remark}

\begin{financialsupport}
  G.~J.~was partially supported by the Swedish Research Council
  (Vetenskapsrådet) Research Project Grant 2022-03748 `Higher structures in
  higher-dimensional homological algebra.' F.~M.~was partially supported by grant
  PID2020-117971GB-C21 funded by MCIN/AEI/10.13039/501100011033.
\end{financialsupport}

\begin{acknowledgements}
  This note is the authors' contribution to the proceedings of the XXI
  International Conference on Representations of Algebras that took place at the
  Shanghai Jiao Tong University in Shanghai, China from July 31 to August 9,
  2024. The authors thank the Scientific Committee for the opportunity of
  presenting their results at the conference, and also thank the Organising
  Committee for their kind hospitality.
\end{acknowledgements}

\printbibliography

\end{document}
